\documentclass[12pt, a4paper]{amsart}

\usepackage[english]{babel}
\usepackage{amssymb,amsmath,amsthm}
\usepackage{verbatim}

\usepackage[all]{xy}
\usepackage{xcolor}

\usepackage{hyperref}
\usepackage{mathrsfs}
\usepackage{color}

\numberwithin{equation}{section}

\newtheorem{dummy}{dummy}[section]

\newtheorem{definition}[dummy]{Definition}
\newtheorem{theorem}[dummy]{Theorem}

\newtheorem{lemma}[dummy]{Lemma}
\newtheorem{proposition}[dummy]{Proposition}

\newtheorem{example}[dummy]{Example}
\newtheorem{question}[dummy]{Question}

\def\A{\mathbb A}
\def\C{\mathbb C}
\def\L{\mathbb L}
\def\P{\mathbb P}

\def\Z{\mathbb Z}

\def\DD{\mathcal D}

\def\HH{\mathcal H}
\def\KK{\mathcal K}
\def\LL{\mathcal L}
\def\MM{\mathcal M}

\def\OO{\mathcal O}

\def\UU{\mathcal U}

\def\LG{\mathrm{LG}}
\def\Db{\DD^{b}}

\def\={\;=\;}
\def\bal{\begin{aligned}}
\def\eal{\end{aligned}}
\def\be{\begin{equation}\label}
\def\ee{\end{equation}}

\def\wt{\widetilde}
\def\ol{\overline}

\DeclareMathOperator{\Gr}{Gr}
\DeclareMathOperator{\Aut}{Aut}
\DeclareMathOperator{\Jac}{Jac}

\DeclareMathOperator{\NS}{NS}
\DeclareMathOperator{\Pic}{Pic}
\DeclareMathOperator{\rk}{rk}

\DeclareMathOperator{\Ker}{Ker}

\title[L-equivalence for dual elliptic quintics]{L-equivalence for degree five elliptic curves,
elliptic fibrations and K3 surfaces}
\author{Evgeny Shinder}
\address{School of Mathematics and Statistics, University of Sheffield, The Hicks Building, Hounsfield Road, Sheffield S3 7RH, UK}
\email{e.shinder@sheffield.ac.uk}
\author{Ziyu Zhang}
\address{Institute for Algebraic Geometry, Leibniz University Hannover, Welfengarten 1, 30167 Hannover, Germany}
\email{zhangzy@math.uni-hannover.de}
\date{}
         


\begin{document}
\maketitle

\begin{abstract}
We construct nontrivial L-equivalence 
between curves of
genus one and degree five, and between elliptic
surfaces of multisection index five.
These results give
the first examples of L-equivalence for curves (necessarily
over non-algebraically closed fields)
and provide a new bit of evidence for the conjectural
relationship
between L-equivalence and derived equivalence.

The proof of the L-equivalence for curves is based on
Kuznetsov's Homological Projective Duality
for $\Gr(2,5)$, and L-equivalence
is extended from genus one curves to elliptic surfaces using
the Ogg--Shafarevich theory of twisting
for elliptic surfaces.

Finally, we apply our results
to K3 surfaces and investigate 
when the two elliptic L-equivalent
K3 surfaces we construct are isomorphic,
using Neron--Severi lattices, moduli spaces of sheaves
and derived equivalence.
The most interesting case is that of elliptic
K3 surfaces of polarization
degree ten and multisection index five,
where the resulting L-equivalence is new.
\end{abstract}

\section{Introduction}

\subsection{The Grothendieck ring of varieties and L-equivalence}

Recall that the Grothendieck ring of varieties $K_0(Var/k)$
is generated as an abelian group by isomorphism 
classes $[X]$ of schemes of finite type $X/k$ modulo the scissor relations
\[
[X] = [U] + [Z]
\]
for every closed $Z \subset X$ with open complement 
$U = X \setminus Z$.
The product structure on $K_0(Var/k)$ is induced by product of schemes. We write $\L \in K_0(Var/k)$
for the class of the affine line $[\A^1]$.

The concept of L-equivalence stems from the recently discovered fact 
that $\L$ is a zero-divisor \cite{Borisov}.
Specifically, for Calabi-Yau threefolds $X$, $Y$
in the so-called
Pfaffian-Grassmannian correspondence, the classes satisfy 
$[X] \ne [Y]$ and
\begin{equation}\label{eq:L-equiv}
\L^n \cdot ([X] - [Y]) = 0,
\end{equation}
where one can take any $n \geqslant 6$ \cite{Borisov, Martin}.
Following \cite{KuznetsovShinder}, we say that smooth projective connected
varieties $X$ and $Y$ are \emph{L-equivalent} if the equation \eqref{eq:L-equiv}
holds for some $n \geqslant 1$,
and we say that $X$ and $Y$ are \emph{nontrivially L-equivalent} if
in addition $[X] \ne [Y]$. If $X$ and $Y$ are not covered by rational
curves and $X$ and $Y$ are not birational then 
an L-equivalence between them is automatically nontrivial
(see e.g. \cite[Proposition 2.2]{KuznetsovShinder}).

There are at least two important reasons why one would want to study L-equivalence.
Firstly, it seems to be closely related to derived equivalence \cite{KuznetsovShinder, IMOU-K3, Kawamata}.
As an evidence for this, the classes of derived categories of L-equivalent
varieties in the Bondal-Larsen-Lunts ring of triangulated categories \cite{BLL}
are equal, and since for Calabi-Yau varieties the derived categories
are indecomposable, it is very likely that nontrivially L-equivalent
Calabi-Yau varieties are actually derived equivalent (see \cite{KuznetsovShinder,
IMOU-K3} for an extended discussion of this relationship).
In fact all currently known examples of pairs of nontrivially L-equivalent varieties are known to be derived equivalent.
These examples include K3 surfaces
\cite{KuznetsovShinder, HassettLai, IMOU-K3, KKM},
Calabi-Yau threefolds
\cite{Borisov, IMOU-G2, BCP},
Calabi-Yau fivefolds \cite{Manivel}
and Hilbert schemes of points on K3 surfaces
\cite{Okawa}.

The second reason to study L-equivalence is the relation to rationality
problems, specifically to that of cubic fourfolds.
Namely, the approach of \cite{GalkinShinder} can be used
to show that very general cubic fourfolds are not rational as soon
as one has sufficient control over the L-equivalence relation.

In this paper we study L-equivalence for genus one curves
and elliptic surfaces, in particular for elliptic K3 surfaces.

\subsection{Genus one curves}

We work over a field of characteristic zero.
Let $X$ be a genus one curve with a line bundle of degree $d$.
For every $k$ coprime to $d$ we can consider the Jacobian $Y = \Jac^k(X)$
which is a fine moduli space
parametrizing degree $k$ line bundles on $X$.
Of course, if $X$ has a rational point, then
all Jacobians $\Jac^k(X)$
are isomorphic to $X$, however
in general this is not the case, and $X$ and $Y$ are typically different torsors
over the same elliptic curve $E = \Jac^0(X)$.

\begin{theorem}\cite{AKW}\label{thm:intro-AKW}
If $k$ and $d$ are coprime, then
genus one curves $X$ and $\Jac^k(X)$ are derived equivalent, and furthermore, every
smooth projective variety $Y$ derived equivalent
to $X$ will be of the form $Y = \Jac^k(X)$
for some $k$ coprime to $d$. 
\end{theorem}

In light of a conjectural relation between L-equivalence
and derived equivalence we may ask the following:

\begin{question}
When are genus one curves $X$ and $Y = \Jac^k(X)$ L-equivalent?
\end{question}

Due to the periodicity relations $\Jac^{k+d}(X) \simeq \Jac^k(X)$,
$\Jac^{-k}(X) \simeq \Jac^k(X)$ and the isomorphism $X \simeq \Jac^1(X)$,
the first nontrivial test case is $d = 5$.
Furthermore, in the $d = 5$ case the only nontrivial coprime Jacobian is 
$Y = \Jac^2(X) \simeq \Jac^3(X)$.
Our first main result is the following:

\begin{theorem}(see Theorem \ref{thm:L-equiv-curves})
If $X$ is a genus one curve with a line bundle of degree $5$ and $Y = \Jac^2(X)$,
then $X$ and $Y$ are L-equivalent, and in general this L-equivalence is nontrivial.
\end{theorem}

More precisely, we show that  \eqref{eq:L-equiv} holds
for $X$ and $Y$ when $n \geqslant 4$ (and does not hold for $n=0$).
This is the first existing construction of nontrivial L-equivalence for curves,
as all the previous constructions were for K3 surfaces or Calabi-Yau varieties
of higher dimension.

As it is often the case with proving L-equivalence we relate the geometry of $X$ and $Y$
to Homological Projective Duality of A.\,Kuznetsov \cite{Kuznetsov}. 
Specifically, as one of the steps in the proof of the theorem above
we prove the following:

\begin{proposition} (see Proposition \ref{prop:deg-5-dual} for the precise statement)
\label{prop:intro-j2}
If $X$ is a genus one curve with a line bundle of degree $5$ and $Y = \Jac^2(X)$, then
$X$ and $Y$ are homologically projectively dual codimension $5$ linear sections of $\Gr(2,5)$.
\end{proposition}

We note the interplay between the moduli space geometry and
the Homological Projective Duality geometry,
in particular
either of the two approaches can be used to show derived equivalence of $X$ and $Y = \Jac^2(X)$.
If one starts with the $Y = \Jac^2(X)$ description, derived equivalence
follows from Theorem \ref{thm:intro-AKW} and if one starts with the Homological Projective Duality
description of Proposition \ref{prop:intro-j2}, derived equivalence
follows from
\cite{Kuznetsov}.

\medskip

To generalize our work and to 
construct L-equivalence of genus one curves in degrees $d > 5$ it seems necessary to study
explicit geometry of the moduli space of curves of genus one and degree $d$.
To describe the geometry of such moduli 
spaces for small $d$
we can use the classical projective models of genus one curves with a degree $d$ divisor:
\begin{itemize}
    \item[] $d = 2$: double covers of $\P^1$ branched in four points
    \item[] $d = 3$: cubic curves in $\P^2$
    \item[] $d = 4$: intersections of two quadrics in $\P^3$
    \item[] $d = 5$: one-dimensional linear sections of a Grassmannian $\Gr(2,5) \subset \P^9$
\end{itemize}

These explicit descriptions show in particular that for $d \leqslant 5$ the 
corresponding moduli spaces are rational.
It is this description in the $d = 5$ case, together with the geometric characterization
of the self-map $\Jac^2$ on the
moduli space of degree $5$ genus one curves,
given in Proposition \ref{prop:intro-j2} 
that allows us to
prove L-equivalence.

We note that the same explicit geometry
of genus one and degree five curves
has been used to study the average size of $5$-Selmer groups
and the average ranks of elliptic curves \cite{BS}.

For $d > 5$ no such explicit description is known, and furthermore it not
known whether the corresponding
moduli spaces are rational or not for large $d$.

\subsection{Elliptic surfaces}

Let $k$ be an algebraically closed field of characteristic zero.
We work with elliptic surfaces without
a section; by a \emph{multisection index} of such a surface
we mean the minimal fiber degree of a multisection.
Our second main result is:

\begin{theorem}(see Theorem \ref{thm:L-equiv-surf})\label{thm:intro-k3}
If $X \to C$ is an elliptic surface of multisection index $5$ and $Y = \Jac^2(X/C)$,
then $X$ and $Y$ are L-equivalent, and in general this L-equivalence is nontrivial.
\end{theorem}

We note that the derived equivalence of $X$ and $Y$ had been proved by Bridegland \cite{Bridgeland}.

We also investigate the case of elliptic K3 surfaces in detail, and answer the question
when the L-equivalence constructed in the Theorem is in fact nontrivial. Here we take $k = \C$.

L-equivalence for K3 surfaces is one of the central open questions in the field.
As a general structural result it is proved by Efimov \cite{Efimov} that every L-equivalence class of K3 surfaces
contains only finitely many isomorphism classes in it.
Previously known cases when nontrivial L-equivalence of derived equivalent K3 surfaces has been constructed are
K3 surfaces of degrees $8$ and $2$ and Picard rank two \cite{KuznetsovShinder}, K3 surfaces
of degree $12$ and Picard rank one \cite{HassettLai, IMOU-K3}, and K3 surfaces of degree $2$ and Picard rank two \cite{KKM}.

Let $X \to \P^1$ be an elliptic K3 surface
of multisection index five,
then $Y = \Jac^2(X/\P^1)$ is also an elliptic K3 surface
(see e.g. \cite[Proposition 11.4.5]{Huy16}),
and by 
Theorem \ref{thm:intro-k3}
these K3 surfaces are L-equivalent.
The next Propositions
explains when $X$ and $Y$ are not isomorphic.

\begin{proposition}(see Proposition \ref{prop:Jac2-K3})
Let $X \to \P^1$ be an elliptic K3 surface of Picard rank two, multisection index $5$ and polarization of degree
$2d$ ($d$ is well-defined modulo $5$), and let $Y = \Jac^2(X/\P^1)$.
    \begin{enumerate}
        \item If $d \equiv 2 \pmod{5}$ or $d \equiv 3 \pmod{5}$, then $X$ and $Y$ are isomorphic.
        \item If $d \equiv 1 \pmod{5}$ or $d \equiv 4 \pmod{5}$, then $X$ and $Y$ are not isomorphic.
        \item If $d \equiv 0 \pmod{5}$, and $X$ is very general in moduli, then $X$
        and $Y$ are not isomorphic.
    \end{enumerate}
\end{proposition}

We note that for every $d$ such K3 surfaces exist and form an $18$-dimensional irreducible
subvariety in the moduli space 
of degree $2d$ polarized
K3 surfaces. Such elliptic K3 surfaces
may have more than one elliptic fibrations (in fact a Picard rank two elliptic K3 has always one or two
elliptic fibrations),
and by an isomorphism of elliptic K3 surfaces
we mean an isomorphism of K3 surfaces, regardless of the elliptic fibration structure. The above Proposition is proved by analyzing lattice theory
of the corresponding K3 surfaces, along the lines of \cite{Stellari, vanGeemen}.

Explicitly, the case (2) of the Proposition covers ellitpic K3 surfaces
of degrees $12$ ($d \equiv 1 \pmod{5})$ and $8$  ($d \equiv 4 \pmod{5})$,
considered previously in \cite{HassettLai, IMOU-K3} and \cite{KuznetsovShinder} respectively.

The K3 surfaces in case (3) can be geometrically described as intersections of
$\Gr(2,5)$, three hyperplanes and a quadric in $\P^9$, and 
containing an elliptic quintic curve
(see Example \ref{ex:K3-deg10}).
This is a genuinely new instance
of nontrivial L-equivalence between K3 surfaces.

\subsection*{Acknowledgements}

We would like to thank Tom Bridgeland, Tom Fisher, Sergey Galkin,
Daniel Huybrechts, 
Alexander Kuznetsov, 
Jayanta Manoharmayum,
C.S.\,Rajan, Matthias Sch\"utt, Constantin Shramov 
and Damiano Testa for helpful discussions and e-mail correspondences.
In particular we thank Alexander Kuznetsov
for explaining to us how to prove the duality in Proposition \ref{prop:deg-5-dual},
and for his comments on a draft of the paper.

We thank the University of Sheffield, Leibniz University Hannover
and the University of Bonn for opportunities for us to travel
and collaborate, and the Max-Planck-Institut f\"ur Mathematik in Bonn
for the excellent
and inspiring
working conditions, where much of this work has been written.

E.S. was partially supported by 
Laboratory of Mirror Symmetry NRU HSE, RF government grant, ag. N~14.641.31.0001.

\section{Dual elliptic quintics}

In this section we work over a field $k$ of characteristic zero.

\subsection{Hyperplane sections of the Grassmannian}

We recall some standard facts about the Grassmannian $\Gr(2,5)$
and its smooth and singular hyperplane sections.

Let $V$ be a five-dimensional vector space; we consider the Pl\"ucker
embedding $\Gr(2,V) \subset \P(\Lambda^2(V)) \simeq \P^9$
and the hyperplane sections $D_\theta := \Gr(2,V) \cap H_\theta$,
parametrized by points of the dual projective space $[\theta] \in \P(\Lambda^2(V^\vee))$,
where $\theta \in \Lambda^2(V^\vee)$ is a nonzero two-form.

By a kernel of a two-form $\theta \in \Lambda^2(V^\vee)$ we mean the subpace
\[
\Ker(\theta) = \{ v \in V: \theta(v \wedge u) = 0
\text{ for all $u \in V$} \}.
\]
For a non-zero form there are two cases:
\begin{enumerate}
\item General case: $\Ker(\theta)$ is one-dimensional.
Then $\theta$ can be written as
$x_{1} \wedge x_2 + x_3 \wedge x_4$ for some basis in $V$.
\item Special case: $\Ker(\theta)$ is three-dimensional. Then $\theta$ is decomposable and can be 
written as $x_1 \wedge x_2$ in some basis.
In other
words $[\theta] \in \Gr(2,V^\vee) \subset \P(\Lambda^2(V^\vee))$.
\end{enumerate}

It is well-known that the two Grassmannians $\Gr(2,V)$ and $\Gr(2,V^\vee)$
are projectively dual in their Pl\"ucker embeddings.
More precisely, we have the following well-known result:

\begin{lemma}\label{lem:dual}
Let $A \subset \Lambda^2(V^\vee)$ be a linear
subspace, and
consider 
its orthogonal subspace
\[
A^\perp = \{ p \in \Lambda^2(V): 
\theta(p) = 0 \text{ for all } \theta \in A\} \subset \Lambda^2(V).
\]
Then $[U] \in \Gr(2,V)$ is a singular point of
$X_A := \Gr(2, V) \cap \P(A^\perp)$ if and only if for every
$\theta \in A$, $\theta(U) = 0$ and for some $\theta_0 \in A$,
$U \subset \Ker(\theta_0)$.

In particular, the hyperplane section $D_\theta$ is singular if and only
if $\theta \in \Gr(2,V^\vee)$, and in this case the singular locus
of $D_\theta$ is isomorphic to $\P^2$.
\end{lemma}
\begin{proof}
The projective tangent space to $\Gr(2,V)$ at a point $[U]$ is 
$\P(U \wedge V) \subset \P(\Lambda^2(V))$, and it follows that
the hyperplane $H_\theta$ is tangent to $\Gr(2,V)$
if and only if $\theta|_{U \wedge V} = 0$, that is $U \subset \Ker(\theta)$.
Thus if $\Ker(\theta)$ is one-dimensional, $D_\theta$ is smooth, and
if $\theta \in \Gr(2,V^\vee)$ so that the $\Ker(\theta)$ is three-dimensional,
$D_\theta$ is singular along $\Gr(2, \Ker(\theta)) \simeq \P^2$.

More generally, if $\theta_1, \dots, \theta_k$ form a basis of $A$,
and $[U] \in X_A = \Gr(2,V) \cap H_{\theta_1} \cap \dots \cap 
H_{\theta_k}$ 
so that all $\theta_i$ vanish on $U$,
then the projective tangent space to $[U]$ at $X_A$ is 
\[
\P(U \wedge V) \cap H_{\theta_1} \cap \dots \cap H_{\theta_k} \subset \P(\Lambda^2(V)),
\]
and this intersection is not transverse if and only if $\theta_1, \dots, \theta_k$
are linearly dependent when restricted to $\P(U \wedge V)$, which
is equivalent to existence of a nonzero form $\theta \in A$
vanishing on $U \wedge V$, or equivalently
$U \subset \Ker(\theta)$.
\end{proof}

\begin{lemma}\label{lem-Gr}
The class in the Grothendieck 
ring of the Grassmannian is 
\[
[\Gr(2,V)] = 1 + \L + 2\L^2 + 2\L^3 + 2 \L^4 + \L^5 + \L^6
\]
and the classes of its smooth and singular hyperplane sections 
$D_\theta = \Gr(2,V) \cap H_\theta$ are given by:
\[\;
D_\theta = 
\left\{\begin{array}{lr}
	   1 + \L + 2\L^2 + 2\L^3 + \L^4 + \L^5, & \text{$\theta \notin
	   \Gr(2,V^\vee)$} \\
	   1 + \L + 2\L^2 + 2 \L^3 + 2\L^4 + \L^5, & \text{
	   $\theta \in \Gr(2,V^\vee)$} \\
       \end{array}
\right.
\]
\end{lemma}
\begin{proof}
The computation for $\Gr(2,5)$ is standard: it is a variety with an affine cell decomposition whose cells are parametrized by Young diagrams fitting into a $3 \times 2$ rectangle, the
codimension of a cell given by the number of blocks in the diagram \cite[Chapter 1.5]{GH}.

We know that a nonzero $2$-form $\theta$ 
on a five-dimensional space has kernel of dimension $1$ (general case) or $3$ (special case) and this distinguishes smooth hyperplane sections from singular ones. Let $K$ be the kernel of $\theta$.

In the smooth case, when $\dim(K) = 1$, 
the two-dimensional subspace $U \subset V$ can either contain $K$ or intersect it trivially; 
thus the subspace $(U + K) / K$ of the four-dimensional space $V/K$ can have dimension $2$ or $1$. 
The space $V / K$ is endowed with a symplectic form $\ol{\theta}$,
and the subspace $(U + K) / K$ is isotropic by construction.
Using the relations in the Grothendieck ring we compute
\[\bal\;
[\Gr(2,5) \cap H_\theta] &= [\P^3] + ([\P^2] - [\P^1])[\LG(2,4)] = \\
&= 1 + \L + 2\L^2 + 2\L^3 + \L^4 + \L^5, 
\eal\]
where we used that the Lagrangian Grassmannian $\LG(2,4)$ is isomorphic
to a three-dimensional split quadric so that $[\LG(2,4)] = 1 + \L + \L^2 + \L^3$ (see e.g. \cite[Example 2.8]{KuznetsovShinder}).

Similarly in the singular case, when $\dim(K) = 3$, 
the two-dimensional subspace $U \subset V$ can either be contained in $K$ or intersect it along a line and considering $(U + K) / K \subset V/K$ 
yields
\[\bal\;
[\Gr(2,5) \cap H] &= [\P^2] + ([\P^3] - [\P^1])[\P^2] = \\
&= 1 + \L + 2\L^2 + 2\L^3 + 2\L^4 + \L^5,
\eal\]
which finishes the proof.
\end{proof}

\begin{proposition}\label{prop:univ-hyper}
For any locally closed subset $S \subset \P(\Lambda^2(V^\vee))$, consider
the universal hyperplane section of $\Gr(2,V)$:
\[
\HH_S := \{ ([U] \in \Gr(2,V), [\theta] \in S): [U] \in H_\theta \} \subset \Gr(2,V) \times S.
\]
Then we have
\[
[\HH_S] = [S] (1 + \L + 2\L^2 + 2\L^3 + \L^4 + \L^5) + 
\L^4 \cdot [S \cap \Gr(2,V^\vee)]. 
\]
\end{proposition}
\begin{proof}
Presenting $S$ as $(S \setminus \Gr(2, V^\vee)) \cup 
(S \cap \Gr(2, V^\vee))$, we see that it suffices to show the
statement when either $S \subset 
\P(\Lambda^2(V^\vee)) \setminus \Gr(2,V^\vee)$
or $S \subset \Gr(2,V^\vee)$.

Let $S \subset \Gr(2,V^\vee)$.
The family of kernels $\Ker(\theta)$, $\theta \in S$ forms a locally-free 
sheaf of rank three over $\HH_S$, and considering the relative position of the fibers
of this sheaf with respect to the fibers of the tautological bundle coming
from $\Gr(2,V)$ allows to repeat the proof of Lemma \ref{lem-Gr}
and to deduce that
\[
[\HH_S] = [S] (1 + \L + 2\L^2 + 2\L^3 + 2\L^4 + \L^5),
\]
which is what we had to prove in this case.

The other case is proved analogously.
\end{proof}

We need one more result regarding incidence rank one sheaves
on hyperplane sections of Grassmannians. Let $V$ be an $n$-dimensional space, 
and let $D \subset \Gr(k,n)$ be the Schubert divisor 
$\sigma_{1,0,\dots,0}$
corresponding to a fixed $(n-k)$-dimensional
linear subspace $W \subset V$, that is
\[
D := \{ [U] \in \Gr(k,n): \dim(U \cap W) \geqslant 1 \} \subset \Gr(k,n).
\]
See \cite[Chapter 1.5]{GH} for the basic properties of the Schubert cycles
$\sigma_{a_1, \dots, a_k}$.

Consider the resolution $\wt{D} \to D$ defined as
\begin{equation}\label{eq:D-res}
\wt{D} := \{ ([U], [l]) \in \Gr(k,n) \times \P(W): 
l \subset U \cap W \}.
\end{equation}

Then $\wt{D}$ is a Grassmannian bundle over
$\P(W)$. We write $h$ for the hyperplane section
on $\P(W)$, as well as for its class on $\wt{D}$,
and we write $H$ for the hyperplane section
on $\Gr(k,V) \subset \P(\Lambda^k(V))$ and its
class on $\wt{D}$.

\begin{lemma}\label{lem:degree}
The $H$-degree of the $c_1(\OO(h)) \in \Pic(\wt{D})$ 
is equal
to the degree of the Schubert cycle $\sigma_{2,0,\dots,0}$ on $\Gr(k,n)$,
that is
\[
c_1(\OO(h)) \cdot H^{n(n-k)-2}
= \sigma_{2,0,\dots,0} \cdot H^{n(n-k)-2}.
\]
\end{lemma}
\begin{proof}
A codimension one linear subspace $W' \subset W$
gives rise to an irreducible
divisor representing $c_1(\OO(h))$:
\[
Z = \{ ([U], [l]) \in \wt{D}: 
l \subset U \cap W' \} \subset \wt{D},
\]
and this divisor maps birationally onto its
image
\[
\{ [U] \in \Gr(k,n): \dim(U \cap W') \geqslant 1 \} \subset D \subset \Gr(k,n).
\]
This subvariety represents the class $\sigma_{2,0,\dots,0}$ in the Chow groups
of the Grassmannian, and it follows that
$H$-degree of $c_1(\OO(h))$ is equal to the $H$-degree
of $\sigma_{2,0,\dots,0}$.
\end{proof}

\subsection{Elliptic quintics, Jacobians and duality}

\begin{definition}
An elliptic quintic is a smooth projective genus one curve which admits a line bundle of degree five.
\end{definition}

By Riemann-Roch theorem a degree five line
bundle $\LL$ on an elliptic quintic $X$ is very ample
and defines an embedding $X \subset \P H^0(X, \LL)^\vee = \P^4$.

\begin{lemma}\label{lem:deg-5}
Let $V$ be a $5$-dimensional $k$-vector space, and $A \subset \Lambda^2(V^\vee)$ be a $5$-dimensional subspace. If $X = \Gr(2,V) \cap \P(A^\perp)$ is a transverse intersection, 
then $X$ is an
elliptic quintic and every elliptic quintic is obtained in this way.
\end{lemma}
\begin{proof}
The first claim follows from the adjunction formula,
while the second one is a classical fact
known as existence of a Pfaffian representation
for an elliptic quintic, see \cite{Fisher}
for a modern exposition.
\end{proof}

For any smooth projective curve $X$ and an integer
$k \in \Z$ we consider the degree $k$ Jacobian $\Jac^k(X)$, defined as the moduli
space of degree $k$ line bundles on $X$. 
If $X$ is an elliptic quintic, then 
by tensoring with the degree $5$ line bundle and by dualizing we obtain the isomorphisms
\[
\Jac^{k+5}(X) \simeq \Jac^k(X), \quad \Jac^{-k}(X) \simeq \Jac^k(X).
\]
Thus in this case 
all Jacobians are isomorphic to one of the
\[
E := \Jac^0(X), \quad X = \Jac^1(X) \simeq \Jac^4(X), \quad Y = \Jac^2(X) \simeq \Jac^3(X).
\]

Here $E$ is an elliptic curve, that is a genus
one curve with a rational point and $X$ and $Y$ are $E$-torsors. 
$E$-torsors are parametrized by the Weil-Chatelet group $H^1(k, E)$
\cite[X.3]{Silverman}. 
If $[X] \in H^1(k,E)$
is the class of the torsor $X$, it is well-known that
for
any $k \in \Z$, $d \cdot [X] = [\Jac^k(X)]$ (see e.g. \cite[Remark 11.5.2]{Huy16}).

In particular,
we see that since $X$ has degree five,
then the order of $[X]$ equals five unless
$X$ has a rational point in which case $[X] = 0$.
Let $Y = \Jac^2(X)$,
then $X \simeq \Jac^2(Y)
\simeq \Jac^3(Y)$.
We call $X$ and $Y$ the \emph{dual elliptic quintics}.
It is clear that if $X$ has a rational point, 
which is always the case when the base field $k$ is algebraically
closed, then $X$ and $Y$ are isomorphic.

We have the following almost converse result.

\begin{lemma}\label{lem:iso}
If $X$ has no rational points and the $j$-invariant
satisfies $j(E) \ne 1728$ then $X$ and $Y$ are not isomorphic. 
\end{lemma}
\begin{proof}
The dual elliptic quintics $X$ and $Y$ give rise to elements $[X], [Y] \in H^1(k, E)$ of order
five, and $[Y] = 2[X]$.

The classes $[X], [Y]$ correspond to
isomorphic genus one curves if and only if $[Y]$ 
lies in 
the $\Aut(E)$-orbit of $[X]$ in $H^1(k, E)$
\cite[Exercise 10.4]{Silverman}.

If we assume that for an automorphism $\sigma \in \Aut(E)$
we have $\sigma([X]) = [Y] = 2[X]$,
the action of $\sigma$ on $H^1(k,E)$ preserves the subgroup $\Z/5$ generated by $[X]$
and we get a surjective 
group homomorphism $\langle \sigma \rangle \to (\Z/5)^\star \simeq \Z/4$.
In particular the order of $\sigma$ should be a multiple of $4$.
On ther other hand since $j(E) \ne 1728$ and $\mathrm{char}(k)=0$,
we have $\Aut(E) = \Z/2$ or $\Aut(E) = \Z/6$,
and no such $\sigma$ exists.

Thus $X$ and $Y$ are not isomorphic.
\end{proof}

We now explain duality between elliptic quintics in terms of projective duality. 

\begin{proposition}\label{prop:deg-5-dual}
Let $V$ be a $5$-dimensional $k$-vector space
and let
$A \subset \Lambda^2(V^\vee)$ be a $5$-dimensional subspace.
We  consider the Grassmannian $\Gr(2,V) \subset \P(\Lambda^2(V))$
and the dual Grassmannian $\Gr(2,V^\vee) \subset \P(\Lambda^2(V^\vee))$. 
For a five-dimensional linear subspace $A \subset \Lambda^2(V^\vee)$ let
\[\bal
X &:= \Gr(2,V) \cap \P(A^\perp)  \\
Y &:= \Gr(2,V^\vee) \cap \P(A).  \\
\eal\]
Assume that $X$ is a smooth transverse intersection,
so that $X$ is a genus one curve.
Then $Y$ is also a smooth transverse intersection and
$X$ and $Y$ are dual elliptic quintics, that
is we have 
\[
Y \simeq \Jac^3(X), \quad 
X \simeq \Jac^2(Y).
\]
\end{proposition}
\begin{proof}
By \cite[Proposition 2.24]{DK} if $X$ is a smooth transverse intersection,
then the same is true for $Y$.

We construct a line bundle
$\MM$ on $X \times Y$.
At each point $([U],[\theta]) \in X \times Y$ we consider the vector space
$\MM_{[U], \theta} := 
U \cap \Ker(\theta)$.
Let us show that this space is one-dimensional.
On the one hand we have $\theta(U) = 0$ so that $U$ can not have trivial
intersection with $\Ker(\theta)$, otherwise
dimension of $\Ker(\theta)$ would be greater than $3$.
 On the other hand $U$ can not be contained in $\Ker(\theta)$, otherwise $[U]$ would be a singular
point of $X$ by Lemma \ref{lem:dual}.

Thus $\MM$, considered as a sheaf
given by the kernel of
\[
p_1^* \left(\UU\big|_X\right) 
\oplus
p_2^* \left(\KK\big|_Y\right)
\to V 
\otimes \OO_{X \times Y}
\]
on $X \times Y$, where $p_1$, $p_2$ are the projections from $X \times Y$
on the two factors,
$\UU \subset V \otimes \OO_{\Gr(2,V)}$ is
the tautological rank two subbundle
on $\Gr(2,V)$ and $\KK \subset V \otimes \OO_{\Gr(2,V^\vee)}$
is the rank three subbundle of kernels
of $2$-forms,
is a locally free sheaf of rank one.

We now compute the bidegree of $\MM$.
For any $\theta \in Y$, 
since $X$ does not intersect the singular
locus of $D_\theta$ (otherwise $X$ would have
been singular), $X$ can be considered
as a curve on the resolution $\wt{D_\theta}$ defined
by \eqref{eq:D-res}.

It follows from definitions that the restriction
$\MM\big|_{X \times \theta}$ is isomorphic
to the restriction of the line
bundle $\OO(-h)$ from $\wt{D_\theta}$ to $X$,
and thus by Lemma \ref{lem:degree} the degree
of $\MM\big|_{X \times \theta}$ is equal up to sign
to the degree of $\sigma_{2,0}$ in $\Gr(2,5)$.
The latter degree is equal to three, as can
be computed using the Pieri formula \cite[Chapter 1.5]{GH}.

The Fourier-Mukai transform defined by $\MM$
is a derived equivalence between $X$ and $Y$
by \cite[Section 4.1, Section 6.1]{Kuznetsov},
which by a standard argument implies that
$X$ and $Y$ are moduli spaces of line bundles on
each other with $\MM$ playing the
role of the universal bundle.

Thus we see that
\[
Y \simeq \Jac^{-3}(X),
\]
and we have $\Jac^{-3}(X) \simeq \Jac^3(X)$ by taking
dual bundles.

Finally, $X \simeq \Jac^2(Y)$ follows by symmetry by repeating the last part of the above argument with the roles of $X$ and $Y$ switched, as the degree of the Schubert cycle
$\sigma_{2,0,0}$ on $\Gr(3,5)$ is equal to two.
\end{proof}

We now deduce L-equivalence of the dual elliptic quintics from their projective duality construction.

\begin{theorem}\label{thm:L-equiv-curves}
Let $X$ and $Y$ be smooth projective dual elliptic quintics. 
Then $X$ and $Y$ are L-equivalent, more precisely we have 
\[
\L^4([X] - [Y]) = 0,
\]
and in general $[X] \ne [Y]$.
\end{theorem}
\begin{proof}
By Lemma \ref{lem:deg-5} 
and Proposition \ref{prop:deg-5-dual}
there exists
a five-dimensional subspace $A \subset \Lambda^2(V^\vee)$ such that 
\[\bal
X &\simeq \Gr(2,V) \cap \P(A^\perp) \\
Y &\simeq \Gr(2,V^\vee) \cap \P(A). \\
\eal\]

We consider the universal hyperplane section $\HH \subset \Gr(2,V) \times \P(A)$:
\[
\HH := \{ U \in \Gr(2,V), \theta \in \P(A)\colon \theta(U) = 0 \} 
\]
and compute its class in the Grothendieck ring
of varieties in two ways.

We apply Proposition \ref{prop:univ-hyper} to
$S := \P(A) \subset \P(\Lambda^2(V^\vee))$ to obtain
\begin{equation}\label{eq-1}
[\HH] = [\P^4] (1 + \L + 2\L^2 + 2\L^3 + \L^4 + \L^5) + \L^4 \cdot [Y]. 
\end{equation}
On the other hand, the morphism $\HH \to \Gr(2,V)$
is Zariski locally-trivial over locally-closed 
subset $\Gr(2,V) \setminus X$ and $X$ with fibers
$\P^3$ and $\P^4$ respectively so that we have
\begin{equation}\label{eq-2}
[\HH] = [\Gr(2,5)][\P^3] + \L^4 \cdot [X]. 
\end{equation}

We compare (\ref{eq-1}) and (\ref{eq-2}).
An easy computation shows that
both $[\P^4] (1 + \L + 2\L^2 + 2\L^3 + \L^4 + \L^5)$
and $[\Gr(2,5)][\P^3]$
are equal to
\[
\L^9 + 2 \L^8 + 4 \L^7 + 6 \L^6 + 7 \L^5 + 7 \L^4 + 6 \L^3 + 4 \L^2 + 2 \L + 1
\]
(for $[\Gr(2,5)]$ see
Lemma \ref{lem-Gr}).
Thus (\ref{eq-1}) and (\ref{eq-2}) together
give
\[
\L^4\cdot([X] - [Y]) = 0.
\]

Finally $X$ and $Y$ are in general not isomorphic
by Lemma \ref{lem:iso}, and since $X$ and $Y$ are not uniruled,
the standard argument shows that $[X] \ne [Y]$ \cite[Proposition 2.2]{KuznetsovShinder}.
\end{proof}

\section{Elliptic surfaces of index five}

In this section $k$ is an algebraically closed field of characteristic zero, and we assume $k = \C$ when
discussing Hodge lattices of K3 surfaces.

\subsection{L-equivalence of elliptic surfaces}

We refer to \cite[Chapter 2]{Dolgachev} for general discussion of elliptic
surfaces and their Jacobians. We recall the basic concepts.
By an \emph{elliptic surface} we mean a smooth
projective surface $X$ with a morphism $\pi: X \to C$
to a smooth projective curve $C$
such that the general fiber of $\pi$
is a genus one curve. We always assume that $X$ is relatively minimal, that is the fibers
of $\pi$ do not contain $(-1)$-curves. 

We do not assume that $\pi$ admits a section.
By the \emph{index} of an elliptic surface we mean the minimal positive degree of 
a multisection of $\pi$.

For every $k \in \Z$ one can consider the relative Jacobian $Y = \Jac^k(X/C)$; $Y$ is another elliptic
surface over the same base curve $C$ defined as the unique 
minimal regular model with 
the generic fiber $\Jac^k(X_{k(C)})$.
As in the genus one curve case, if $X$ admits a section, then all Jacobians $\Jac^k(X/C)$ are isomorphic to $X$ over $C$. 

\begin{lemma}
\label{lem:isom-fibers}
If $X \to C$ is an elliptic surface
and $Y = \Jac^k(X/C)$, then
for every point $c \in C$, the reduced fibers
$(X_c)_{red}$ and $(Y_c)_{red}$ are isomorphic.
\end{lemma}
\begin{proof}
This follows from \cite[Chapter 2, Proposition 1 and 2]{Dolgachev}.
\end{proof}

We now consider the case when the multisection
index of an elliptic surface
$X \to C$ is equal to five, 
and analogously to the genus one curve case
we call $X$ and $Y = \Jac^2(X/C) 
\simeq \Jac^3(X/C)$ the dual elliptic fibrations.

\begin{theorem}\label{thm:L-equiv-surf}
Let $X \to C$ be an elliptic fibration of index five over an algebraically closed
field of characteristic zero, and let $Y = \Jac^2(X/C)$.
Then $X$ and $Y$ are $L$-equivalent, more precisely we have
\[
\L^4([X] - [Y]) = 0.
\]
\end{theorem}
\begin{proof}
    Let $X_{k(C)}$, $Y_{k(C)}$ be the generic fibers of $X$ and $Y$.
    By Theorem \ref{thm:L-equiv-curves}, we have $\L^4([X_{k(C)}] - [Y_{k(C)}]) = 0$ in $K_0(Var/k(C))$. Therefore by \cite[Proposition 3.4]{NicaiseSebag}
    there exists a non-empty open set $U \subset C$ such that 
    \[
    \L^4([X_U] - [Y_U]) = 0,
    \]
    in $K_0(Var/k)$, where $X_U$, $Y_U$ are preimages of $U$ in $X$ and $Y$ respectively.
    Let $C \setminus U = \{ c_1, \dots, c_n \}$, 
    then $X_{c_i}$ and $Y_{c_i}$ are isomorphic for each $i$ 
    by Lemma \ref{lem:isom-fibers}.
    In particular $\L^4([X_{c_i}]-[Y_{c_i}]) = 0$; summing everything together we obtain the desired L-equivalence statement.
\end{proof}

In the next section we show that elliptic K3 surfaces of index five and Picard
rank two provide examples when $X$ and $Y$
are not isomorphic, see Proposition \ref{prop:Jac2-K3},
so that $[X] \ne [Y]$ (see e.g. \cite[Proposition 2.8]{KuznetsovShinder}).

\subsection{Elliptic K3 surfaces of Picard rank two}

We consider elliptic K3 surfaces over
$k = \C$.
Recall that for a K3 surface
$\NS(X) \simeq \Pic(X)$ is a free
finitely generated abelian group whose
rank is called the Picard rank of $X$.
Intersection pairing gives
$\NS(X)$ a structure of a lattice.
See \cite[Chapter 14]{Huy16}
for an introduction to lattices.
We write $U$ for the hyperbolic plane, and 
$N(X)$ for the extended Neron--Severi lattice $N(X) = U \oplus \NS(X)$ under the Mukai pairing.
We say that two indefinite lattices have the same genus
if they have the same rank, signature and discriminant groups.

We only consider projective K3 surfaces,
that is the ones admitting a polarization.
We think of polarization as a class of an ample
divisor in $\NS(X)$.
Since by degree reasons the class of a polarization is linearly independent to
the class of the fiber of an elliptic fibration, the minimal Picard rank of an elliptic
K3 surface is equal to two.
Good references about such K3 surfaces
are papers of Stellari \cite{Stellari} and van Geemen \cite{vanGeemen},
and \cite[Chapter 11]{Huy16}.

\begin{lemma}\label{lem:ell-K3-dt}\cite[Remark 4.2]{vanGeemen}
Let $X$ be an elliptic K3 surface of index $t > 0$ and of Picard rank two.
Let $F \in \NS(X)$ be the class of the fiber.
Then there exists a polarization $H$ such that $H \cdot F = t$, 
and $H$, $F$ form a basis of $\NS(X)$.
\end{lemma}
\begin{proof}
Let us show that $F \in \NS(X)$ is a primitive class. Indeed, if $F = mC$,
for $m \ge 1$, then $C$ will be an effective
divisor contained in a fiber.
Since we assume that Picard rank of
$X$ is two, all fibers are irreducible,
and $m = 1$.

Since $F$ is a primitive class, there exists
$D \in \NS(X)$ such that $D, F$ form a basis of $\NS(X)$. Up to replacing $D$ by $-D$
we may assume that $D \cdot F = t$. 
A simple computation shows that the only possible $(-2)$-classes in $\NS(X)$ are given by $\pm (D + \frac{2-D^2}{2t}F)$, hence there is at most one $(-2)$-curve in $X$.

We consider $H = D + nF$. It is clear that
	$$ H^2 = D^2 + 2n t >0 $$
	for $n \gg 0$. If $C$ is a $(-2)$-curve, then
	$$ H \cdot C = D \cdot C + n F \cdot C > 0 $$
	for $n \gg 0$ since $C$ is not in any fiber (otherwise the Picard rank of $X$ would be at least
	three). 	Hence $H$ is ample for $n \gg 0$ by \cite[Proposition 2.1.4]{Huy16}. 
\end{proof}

For a pair of integers $t > 0$
and $d \in \Z$
we consider a rank two lattice $\Lambda_{t,d}$
with basis $H$, $F$ and pairing defined by
\begin{equation}\label{eq:lattice}
\begin{pmatrix}
		2d & t \\
		t & 0. \\
\end{pmatrix}
\end{equation}

\medskip

There always exist projective
K3 surfaces with $\NS(X) \simeq \Lambda_{t,d}$
\cite[Corollary 14.3.1]{Huy16}.
Any such K3 surface is elliptic because
$\NS(X)$ contains a square-zero class \cite[Proposition 11.1.3]{Huy16}.
Furthermore since
the embedding of $\Lambda_{t,d}$ into a K3 lattice
is unique up to isomorphism by \cite[Corollary 14.3.1]{Huy16} 
the locus of these K3 surfaces
is an irreducible locally closed subset
of dimension $18$ in the moduli space of all
degree $2d$ polarized K3 surfaces. 

Note that $t$ is a well-defined invariant
of $\Lambda_{d,t}$, as the discriminant
of \eqref{eq:lattice} is $-t^2$.
The following result describes
the complete set
of invariants of $\Lambda_{d,t}$
in the case when $t$ is an odd prime.

\begin{proposition}[van Geemen, Stellari]
\label{prop:Lambda-dt}
Let $t > 0$ be an odd prime, and let
$d, d' \in \Z$.

\noindent(1) $\Lambda_{d,t}$ is isomorphic
to $\Lambda_{d',t}$ if and only if 
$d \equiv d' \pmod{t}$
or $d d' \equiv 1 \pmod{t}$.

\noindent(2) 
$O(\Lambda_{d,t}) = \{\pm 1\}$
if $d \not\equiv \pm 1 
\pmod{t}$
and $O(\Lambda_{d,t}) = \Z/2 \times \Z/2 =
\{\pm 1, \pm J\}$
where $J$ is the isometry swapping
the two isotropic classes if
$d \equiv \pm 1 
\pmod{t}$.

\noindent(3) The discriminant group
$A_{d,t} = \Lambda_{d,t}^\star / \Lambda_{d,t}$ is $(\Z/t)^2$
if $t$ divides $d$, 
and for $\gcd(d,t) = 1$,
it is $A_{d,t} = \Z/t^2$ with
the square of the generator given
by $\frac{-2d}{t^2}$.

\noindent (4) $\Lambda_{d,t}$, $\Lambda_{d',t}$
are in the same genus if and only 
$d' \equiv k^2 d \pmod{t}$ for some 
integer $k$ coprime to $t$.
\end{proposition}
\begin{proof}
(1) is \cite[Proposition 3.7]{vanGeemen}.
and (2) is \cite[Lemma 4.6]{vanGeemen}.
The result in (3) is easy for $t | d$ 
as we can assume $d = 0$.
For $\gcd(d,t) = 1$, 
(3) is the computation in the proof 
of \cite[Lemma 3.2 (ii)]{Stellari}.
(4) is \cite[Lemma 3.2 (ii)]{Stellari}. 
\end{proof}

\begin{example}\label{ex:Lambda-5}
If $t = 5$, then there are four isomorphism classes of lattices $\Lambda_{5,d}$:
\[
\Lambda_{5,0}, \Lambda_{5,1}, \Lambda_{5,2} \simeq \Lambda_{5,3}, \Lambda_{5,4}.
\]

The discriminant group
$A_{d,t} = \Lambda_{d,t}^\star / \Lambda_{d,t}$
for $\Lambda_{5,0}$ is $\Z/5 \oplus \Z/5$,
and it is $\Z/25$ in the other cases.

The lattices $\Lambda_{5,1}$ and $\Lambda_{5,4}$ are in the same genus, whereas the other lattices have only one isomorphism class in each genus.

Finally, the lattices $\Lambda_{1,5}$, $\Lambda_{4,5}$, $\Lambda_{5,5}$ admit an isometry $J$ permuting the 
two isotropic classes, and the isometry group is $\Z/2 \times \Z/2 = \{\pm 1\} \times \{\pm J\}$, whereas the lattice $\Lambda_{2,5} \simeq \Lambda_{3,5}$ has the isometry group 
$\Z/2 = \{\pm 1\}$.
\end{example}

Explicitly one can get a K3 surface with $\NS(X) = \Lambda_{t,d}$
by taking a general K3 surface containing a degree $t$ elliptic curve.

\begin{example}
\label{ex:K3-8}
A very general degree $8$ K3 surface
$X \subset \P^5$ which contains a normal rational
curve $C$ of degree three, has $H^2 = 8$, $C \cdot H = 3$, $C^2 = -2$, so that
\[
\NS(X) \simeq
\begin{pmatrix}
8 & 3 \\
3 & -2
\end{pmatrix}.
\]
Such a K3 surface admits an elliptic fibration
provided by the pencil $F = H - C$, which consists
of the residual elliptic quintics in the hyperplane
sections of $X$ through $C$ and it is easy to compute
that we have
\[
\NS(X) \simeq \Lambda_{5,4}.
\]
We note that $X$ admits a unique elliptic fibration \cite[4.7]{vanGeemen}.
\end{example}

\begin{example}\label{ex:K3-deg10}
A general degree $10$ K3 surface $X$ is a complete intersection of a Grassmannian
$\Gr(2,5) \subset \P^9$ with three hyperplanes and a quadric \cite[Corollary 0.3]{Mukai-FanoK3}.

As soon as $X$ contains a normal elliptic quintic curve $F \subset \P^4$, it will admit an elliptic
fibration of index five, and generically we have
\[
\NS(X) \simeq
\begin{pmatrix}
10 & 5 \\
5 & 0
\end{pmatrix}.
\]
in the basis $H$, $F$. In fact if we write $F' = H - F$, we see that $\NS(X)$ is isomorphic to $\Lambda_{5,0}$.
We note that $F'$ gives rise to a second elliptic fibration structure on $X$, cf. \cite[4.7]{vanGeemen}.
\end{example}

We prepare to address the question when $\Jac^k(X/\P^1)$ 
and $X$ are isomorphic.

\begin{lemma}\label{lem:NS-Jk}
If $X$ is a K3 surface with $\NS(X) = \Lambda_{t,d}$, 
and $\gcd(t,k) = 1$, then
$\NS(\Jac^k(X/\P^1)) = \Lambda_{t,d \cdot k^2}$
for any elliptic fibration on $X$. 
\end{lemma}
\begin{proof}
Let $N(X) = U \oplus \NS(X)$ be the extended Neron-Severi lattice
and let $e_1, e_2$ be a basis of $U$ 
consisting of two isotropic vectors with $e_1 \cdot e_2 = -1$.
Then $v = F + k e_2 \in N(X)$ is the Mukai vector
giving rise to the moduli space $Y = \Jac^k(X/\P^1)$ \cite[Example 16.2.4]{Huy16}.
Using \cite[Theorem 1.4]{Mukai-K3} we have
\[
\NS(Y) = v^\perp / v.
\]
Explicitly we have 
\[
v^\perp = \langle F, e_2, kH + t e_1 \rangle
= \langle v, e_2, kH + t e_1 \rangle,
\]
so that
\[
v^\perp / v = \langle e_2, kH + te_1 \rangle,
\]
and the intersection form on this lattice
is isomorpic to $\Lambda_{t,d \cdot k^2}$.
\end{proof}

We need the following result,
which describes the group of Hodge
isometries of the transcendental
lattice for a sufficiently general
K3 surface. This group is important
for studying derived equivalence
between K3 surfaces.
In particular it 
appears in the
counting formula for the
number of Fourier-Mukai partners 
\cite[Theorem 2.3]{HLOY}.
In the proof we follow the strategy of
\cite[Proposition B.1]{Oguiso}.

\begin{lemma}\label{lem:G-triv}
If $X$ has Picard rank $\rho < 20$ and $X$
is very general in the moduli space of K3 surfaces polarized by a fixed sublattice $\NS(X)$ of the K3 lattice, then the group
of Hodge isometries of the transcendental lattice $T_X$ is $\{\pm 1\}$.
\end{lemma}
\begin{proof}
By the Torelli theorem for K3 surfaces, a (marked) K3 surface polarized by $\NS(X)$ is determined by a holomorphic $2$-form $\sigma_X \in T_X\otimes\mathbb{C}$, considered up to scalar. Since the choice of the form $\sigma$ is given by the condition 
\begin{equation}
    \label{eqn:choice_sigma}
    \sigma_X^2 = 0 \quad \text{ and } \quad \sigma_X\overline{\sigma}_X>0,
\end{equation}
and for a very general choice of $\sigma_X$ satisfying \eqref{eqn:choice_sigma}, $\sigma_X^\perp$ in $T_X\otimes\mathbb{C}$ contains no non-trivial integral class, we conclude that the moduli of (marked) K3 surface polarized by $\NS(X)$ has
dimension $\rk(T_X) - 2$.

Let us fix an isometry $g$ of $T_X$, and
assume that $g$ induces a \emph{Hodge} isometry of $T_X$ for the K3 surface
$X$ corresponding to $\sigma_X$.
We use \cite[Proposition B.1]{HLOY}.
For any choice of $\sigma_X$, 
the group of Hodge isometries of $T_X$ is a finite cyclic group of even order $2m$,
and without loss of generality we may assume that $g$ is a generator of this group. 
Furthermore in this case
$g$
acts on $\sigma_X$ via multiplication by a primitive $2m$-th root of unity. 
Finally, $T_X \otimes \mathbb{C}$ decomposes into a direct sum of eigenspaces of $g$ as
\begin{equation}
    \label{eqn:decomp_Tx}
    T_X \otimes \mathbb{C} = \bigoplus_{\xi} V_\xi,
\end{equation}
where $\xi$ runs over all primitive $2m$-th roots of unity, and the dimension of each eigenspace $V_\xi$ is $\rk{T_X}/\varphi(2m)$ with $\varphi(-)$ being the Euler function
(see \cite[Steps 4, 5 in the proof of Proposition B.1]{HLOY}). 
Since $\sigma_X$ is an eigenvector for $g$, we have $\sigma_X \in V_\xi$ for some $\xi$. It follows that the moduli of such K3 surfaces has dimension at most $\rk(T_X)/\varphi(2m) -1$.

By assumption $\rho < 20$, so that we have $\rk(T_X)>2$. If $m>1$, then $\varphi(2m) \geqslant 2$ and
$$ \rk(T_X)/\varphi(2m) -1 \leqslant \rk(T_X)/2 -1 < \rk(T_X) - 2, $$ 
where the right-hand-side is the dimension of the moduli space of K3 surfaces
polarized by $\NS(X)$
and the left-hand-side is the dimension of the closed subvariety in the moduli
where $g$ becomes the generator for the group of Hodge isometries.
This means that unless $g = \pm 1$, $g$ is not a Hodge isometry of 
$T_X$ of a general K3 surface in the moduli.

Since the group of isometries of $T_X$
is countable, very general choices of $\sigma_X$ would
give K3 surfaces with the group of Hodge isometries of $T_X$ equal to $\{\pm 1\}$.
\end{proof}

We now consider the multisection index $5$ case. According to Lemma \ref{lem:ell-K3-dt},
an elliptic K3 surface with Picard rank two will have Neron--Severi lattice
isomorphic to one of the $\Lambda_{5,d}$, where $d$ is considered modulo $5$
See Example \ref{ex:Lambda-5} for more details about these lattices.

\begin{proposition}\label{prop:Jac2-K3}
    Let $X$ be an elliptic K3 surface with $\NS(X) = \Lambda_{5,d}$, and let $Y = \Jac^2(X/\P^1)$. 
    \begin{enumerate}
        \item If $d = 2$ or $d = 3$, then $X$ and $Y$ are isomorphic.
        \item If $d = 1$ or $d = 4$, then $X$ and $Y$ are not isomorphic.
        \item If $d = 0$, and $X$ is very general in moduli, then $X$
        and $Y$ are not isomorphic.
    \end{enumerate}
\end{proposition}
\begin{proof}
\noindent (1) It suffices to show that $X$ does not have nontrivial Fourier-Mukai partners. We note that
by Proposition \ref{prop:Lambda-dt} (1) and (4), $\Lambda_{5,2} \simeq \Lambda_{5,3}$
is the only isometry class of a lattice in its genus.
Hence the counting formula for Fourier-Mukai partners \cite[Theorem 2.3]{HLOY} has
only one term and since by Proposition \ref{prop:Lambda-dt} (2) 
the orthogonal group
$O(\Lambda_{5,2})$ consists of $\pm 1$,
this term is equal to one. 

\noindent (2) By Lemma \ref{lem:NS-Jk}, taking
$\Jac^2$ interchanges the Neron--Severi lattices $\Lambda_{5,1}$
and $\Lambda_{5,4}$,
and since these lattices
are not isomorphic, $X$ and $Y$ are not isomorphic.

\noindent (3) 
If $X$ and $Y$ are isomorphic, then the Fourier-Mukai transform $\Phi: \Db(X) \simeq \Db(\Jac^2(X/\P^1))$
corresponding to the moduli space $\Jac^2(X/C)$ on $X$ 
induces a Hodge isometry of 
$H^\star(X, \Z)$ taking one Mukai vector to the other \cite[Section 16.3]{Huy16}.

Consider the extended Neron--Severi lattice $N(X) = U \oplus NS(X)$,
where we choose a basis $e_1, e_2$ 
for $U$ consisitng of two isotropic vectors satisfying $e_1 \cdot e_2 = -1$.
The action of $\Phi$
takes $e_1$ 
(Mukai vector for moduli space $X$ on $X$) to $F + 2e_2$
(Mukai vector
for moduli space $Y$ on $X$).

We note that one such isometry $g_0 \in O(N(X))$ is
\begin{equation}\label{eq:g0-def}
\bal
e_1 &\mapsto 2e_2 + F \\
e_2 &\mapsto -2e_1 - H \\
H &\mapsto 2H + 5e_1 \\
F &\mapsto -2F - 5e_2 \\
\eal
\end{equation}
and any other isometry $g$ mapping $e_1$ to $F + 2e_2$
will have the form
\[
g = g_0 \cdot h,
\]
where $h \in O(N(X), e_1)$ is an isometry of $N(X)$ fixing $e_1$.

We now consider the action of $g$ on the discriminant group $N(X)^\star / N(X) \simeq
NS(X)^\star / NS(X) \simeq A_{5,0} = \Z/5 \oplus \Z/5$ generated
by $\frac15 H, \frac15 F$ (cf. Proposition \ref{prop:Lambda-dt} (3)).
Since we assume that the action of $g$ is induced by a Hodge isometry
of $H^\star(X, \Z)$, the action of $g$ on the discriminant group is the same
as the action induced by a Hodge isometry of $T_X$.
By Lemma \ref{lem:G-triv} for general $X$ this
action on the discriminant group is $\pm 1$.

We note that the action of $O(N(X), e_1)$ on the discriminat group
factors through $O(e_1^\perp / e_1) = O(NS(X))$,
so by \cite[Lemma 4.6]{vanGeemen} its action is given by one of the matrices
\[
\begin{pmatrix}
1 & 0 \\
0 & 1 \\
\end{pmatrix}, 
\begin{pmatrix}
-1 & 0 \\
0 & -1 \\
\end{pmatrix}, 
\begin{pmatrix}
0 & 1 \\
1 & 0 \\
\end{pmatrix}
\begin{pmatrix}
0 & -1 \\
-1 & 0 \\
\end{pmatrix}.
\]

On the other hand we see from \eqref{eq:g0-def} that the action of $\pm g_0$ on $A_{5,0}$
does not belong to the subgroup above. Therefore there is no element $g \in O(N(X))$
which maps $e_1$ to $F + 2e_2$ and is induced by a Hodge isometry of $H^\star(X,\Z)$.
\end{proof}


\providecommand{\arxiv}[1]{{\tt{arXiv:#1}}}

\medskip
\medskip


\end{document}